\documentclass{amsart}
\usepackage{amssymb}

\usepackage{color}
\usepackage[colorlinks]{hyperref}
\usepackage{pb-diagram,lamsarrow,pb-lams}

\theoremstyle{plain}
\newtheorem{thm}{Theorem}[section]

\newtheorem{cor}[thm]{Corollary}
\newtheorem{lem}[thm]{Lemma}
\newtheorem{prop}[thm]{Proposition}

\theoremstyle{definition}
\newtheorem*{defn}{Definition}

\theoremstyle{remark}
\newtheorem{rem}[thm]{Remark}

\newcommand{\thmref}[1]{Theorem~\ref{#1}}

\newcommand{\lemref}[1]{Lemma~\ref{#1}}

\newcommand{\propref}[1]{Proposition~\ref{#1}}

\newcommand{\corref}[1]{Corollary~\ref{#1}}

\def\RR{{\mathbb R}}

\def\BBB{{\mathcal B}}

\def\PPP{{\mathcal P}}

\def\kkkk{{\mathfrak K}}

\newcommand{\booktitle}[1]{\textit{#1}}

\newcommand{\setof}[1]{\left\{#1\right\}}

\def\nabla{\triangledown}

\def\+{\oplus}

\newenvironment{enroman}[0]{

\begin{enumerate}
}{
\end{enumerate}

}
\newcounter{ovaj}
\newcommand{\navedi}[1]{\setcounter{ovaj}{#1}({\bf\roman{ovaj}})}

\begin{document}


\title{Homotopy characterization of ANR mapping spaces}
\author{Jaka Smrekar}
\address{Fakulteta za matematiko in fiziko, Jadranska ulica 19, SI-1111 Ljubljana, Slovenia}
\email{jaka.smrekar@fmf.uni-lj.si}

\thanks{The author was supported in part by the M\v{S}Z\v{S} of the Republic of Slovenia research program
No.~P1-0292-0101-04 and research project No.~J1-6128-0101-04.}
\subjclass[2000]{Primary 55P99}
\date{August 29, 2007}
\keywords{ANR, function space, compact open topology, homotopy type of a CW complex, P-embedding.}


\begin{abstract}
Let $Y$ be an absolute neighbourhood retract (ANR) for the class of metric spaces and let $X$
be a Hausdorff space. Let $Y^X$ denote the space of continuous maps from $X$ to $Y$ equipped
with the compact open topology. It is shown that if $X$ is a CW complex then $Y^X$ is an ANR for the
class of metric spaces if and only if $Y^X$ is metrizable and has the homotopy type of a CW complex.
The same holds also when $X$ is a compactly generated hemicompact space (metrizability assumption is
void in this case).
\end{abstract}

\maketitle


\section{Introduction}

Let $X$ and $Y$ be two spaces having the homotopy type of a CW complex.
Let $Y^X$ denote the space of continuous maps $X\to Y$ equipped with the compact open topology.
The author of this note has extensively investigated the question of when also $Y^X$
has the homotopy type of a CW complex (see \cite{smrekar} and \cite{smr_CRM}).

Shifting the viewpoint to absolute neighbourhood retracts (ANRs) for metric spaces,
one can ask the following question: If $X$ is a CW complex and $Y$ is an ANR,
when is $Y^X$ an ANR? (Note that if $X$ is uncountable, $Y^X$ need not even be metrizable.)

A topological space has the homotopy type of a CW complex if and only if
it has the homotopy type of an ANR (see Milnor \cite{milnor}). However, it is not difficult to find
examples of spaces that are not ANRs but have the homotopy type of a CW complex.

It turns out that if $X$ is any CW complex and $Y$ is an ANR, then
the space $Y^X$ is an ANR whenever it is metrizable and has the homotopy
type of a CW complex:

\begin{thm}\label{main}
Let $X$ be a CW complex and let $Y$ be an ANR for metric spaces.
The following are equivalent.
\begin{enroman}
	\item	The space $Y^X$ is an ANR for metric spaces.
	\item	The space $Y^X$ is metrizable and has the homotopy type of a CW complex.
	\item	The space $Y^X$ is metrizable and semilocally contractible.
\end{enroman}
\end{thm}

(Recall that a space $Z$ is {\it semilocally contractible} if each point
has a neighbourhood that is contractible within $Z$.)

A Hausdorff space $X$ is called {\it hemicompact} if $X$ is the union
of countably many of its compact subsets $\setof{K_i\,\vert\,i}$ which dominate
all compact subsets in $X$. This means that for each compact $K\subset X$ there
exists $i$ with $K\subset K_i$.

Observe that a CW complex is hemicompact if and only if it is countable.
However, a countable CW complex has a distinguished class of compact subsets,
namely, finite subcomplexes. Those are always cofibered in the total space.
No such thing is required for a general hemicompact space.

A space $X$ is {\it compactly generated} if the compact subspaces determine its topology.
For compactly generated hemicompact domain spaces the following is true:

\begin{thm}\label{submain}
Let $X$ be a compactly generated hemicompact space and let $Y$ be an ANR for metric spaces.
Then the following are equivalent.
\begin{enroman}
	\item	The space $Y^X$ is an ANR for metric spaces.
	\item	The space $Y^X$ has the homotopy type of a CW complex.
\end{enroman}
\end{thm}

The proofs of \thmref{main} and \thmref{submain} are quite different. The first
uses the Cauty-Geoghegan's characterization of ANRs (see Cauty \cite{cauty_ANR}),
the second uses Morita's homotopy extension theorem for $P_0$-embeddings (see Morita \cite{morita}).
It is not possible to use the proof of \thmref{main} for \thmref{submain},
and it seems difficult to use the proof of \thmref{submain} for \thmref{main}
in case $X$ is an uncountable CW complex.

{\bf Acknowledgement.}
The author is indebted to Atsushi Yamashita for kindly sending him the preprint \cite{yamashita}
where the question of equivalence of \navedi{1} and \navedi{2} of \thmref{main} (for countable
domains) was posed implicitly. Inspiration for existence of \thmref{submain} was also found there.
Moreover, the preprint shows that the two theorems can be useful in considering function spaces
which are Hilbert manifolds.


\section{Proof of \thmref{main}}

For subsets $A$ of the domain space and $B$ of the target space we let $G(A,B)$ denote
the set of all maps $f$ that map the set $A$ into the set $B$. For topological spaces
$X$ and $Y$, we can take as subbasis of the compact open topology on $Y^X$ the collection
$\PPP$ of all $G(K,V)\subset Y^X$ with $K$ a compact subset of $X$ and $V$ an open subset of $Y$.

We will employ the characterization of absolute neighbourhood retracts conjectured by Ross Geoghegan
and proven by Robert Cauty:

\begin{thm}[Cauty \cite{cauty_ANR}, `Th\'eor\`eme'] \label{th_cauty}
A metrizable space $Z$ is an ANR for metric spaces if and only if each open subset of $Z$
has CW homotopy type.\qed
\end{thm}

\begin{proof}[Proof of Theorem 1.1]
The equivalence of \navedi{2} and \navedi{3} follows from \cite{smr_CRM}, Theorem 2.2.1.
Evidently \navedi{1} implies \navedi{2}. 

Assume that $Y^X$ has the homotopy type of a CW complex. We claim that each
open subset $U$ of $Y^X$ has the homotopy type of a CW complex.

Let $K_1,\dots,K_n$ be compact subsets of $X$ and let $V_1,\dots,V_n$ be open subsets of $Y$.
We show first that the finite intersection $G(K_1,V_1)\cap\dots\cap G(K_n,V_n)$ has the homotopy type of a CW complex.
To this end, let $L$ be a finite subcomplex of $X$ containing the compact union $K_1\cup\dots\cup K_n$.

The restriction mapping $R\colon Y^X\to Y^L$ assigning to each $f\colon X\to Y$ its restriction $f\vert_L\colon L\to Y$
is a Hurewicz fibration because $L\hookrightarrow X$ is a cofibration.
Denote $\widetilde W=G(K_1,V_1)\cap\dots\cap G(K_n,V_n)\subset Y^X$ and
$W=G(K_1,V_1)\cap\dots\cap G(K_n,V_n)\subset Y^L$. The preimage of $W$ under $R$ is exactly $\widetilde W$
hence also $r=R\vert_{\widetilde W}\colon \widetilde W\to W$ is a Hurewicz fibration. Since by a result of Kuratowski
(see also Marde\v{s}i\'c and Segal \cite{mardesic-segal}, Theorem I.3.4), the space $Y^L$ is an ANR, the open
subset $W$ of $Y^L$ has the homotopy type of a CW complex by \thmref{th_cauty}. We claim that all fibres of $r$ have CW
homotopy type as well.

To this end, note that for each $\varphi\in W$, the fibre $F_\varphi$ of $r$ over $\varphi$ in $\widetilde W$
is precisely the fibre of $R$ over $\varphi$ in $Y^X$. As noted above, $Y^L$ is an ANR and hence has CW homotopy
type. By assumption, $Y^X$ has CW homotopy type. By Stasheff \cite{stasheff}, Corollary (13), $F_\varphi$ has
the homotopy type of a CW complex. Thus since the base space $W$ and all fibres of $r$ have CW homotopy type,
also $\widetilde W$ has CW homotopy type by Stasheff \cite{stasheff}, Proposition (0) (see also Sch\"on \cite{schon}).

Let $\BBB$ denote the collection of sets $G(K_1,V_1)\cap\dots\cap G(K_n,V_n)\subset Y^X$
for all possible choices of $n$, $K_i$ (compact), and $V_i$ (open). This is to say that $\BBB$ is the standard
basis for the topology on $Y^X$ associated to the subbasis $\PPP$. Note that $\BBB$ is closed under
formation of finite intersections and that each member of $\BBB$ has CW homotopy type by the above.

Let $U$ be an arbitrary open subset of $Y^X$ and let $\BBB_U$ be the set of those elements of $\BBB$ that
are contained in $U$. Being metrizable, the space $Y^X$ is hereditarily paracompact, hence $\BBB_U$ is
a numerable open covering of the space $U$. By tom Dieck \cite{tdieck}, Theorem 4, it follows that $U$
has the homotopy type of a CW complex, as claimed. An application of \thmref{th_cauty} completes the proof.
\end{proof}

\begin{rem}
The proof leans on the fact that in a CW complex $X$, each compact set $K$ is contained in another
compact subset $L$ for which the inclusion $L\hookrightarrow X$ is a closed cofibration. Moreover,
the topology on $X$ is determined by its compact subsets. The proof of \thmref{main} generalizes
trivially to domain spaces $X$ with these two properties.
\end{rem}

By virtue of \thmref{main}, one can find in \cite{smr_CRM} a number of results implying that 
certain function spaces either are or aren't ANRs. In addition, the results there show that
the problem of determining whether $Y^X$ has CW homotopy type is very hard.

\begin{cor}\label{cmilnor}
Let $X$ be a connected countable CW complex.
Then $X$ is homotopy dominated by a finite CW complex if and only if $\pi_1(X)$ is finitely presentable
and $Y^X$ is an ANR for all ANR spaces $Y$.
\end{cor}

\begin{proof}
Follows from \thmref{main} and Theorem 4.5.3 of \cite{smr_CRM}. Note that the metrizability condition
is void here.
\end{proof}

\begin{rem}
Under additional restrictions on $X$, \corref{cmilnor} was obtained independently by Yamashita in \cite{yamashita}.
\end{rem}


\section{Proof of \thmref{submain}}

To prove \thmref{submain} we use the fact that ANRs for metric spaces are precisely the
absolute neighbourhood extensors for metric spaces. Given the hypotheses of \navedi{2}
of \thmref{submain}, therefore, we need to show that for every pair $(Z,A)$ with $Z$
metric and $A$ closed in $Z$, every continuous function $f\colon A\to Y^X$ extends
continuously over a neighbourhood of $A$ in $Z$. The idea of the proof is very simple, and
we outline it first in case $X$ is, in addition, locally compact. The technical
details for the general case will follow below.

{\bf Outline of proof.} Since $Y^X$ has the homotopy type of an ANR, a continuous map $f\colon A\to Y^X$ always
admits a neighbourhood extension up to homotopy. That is, there exist a continuous map $g\colon U\to Y^X$
where $U$ is open and contains $A$, and a homotopy $h\colon A\times[0,1]\to Y^X$ beginning in $g\vert_A$
and ending in $f$. The maps $g$ and $h$, respectively, induce continuous adjoints $\Hat g\colon U\times X\to Y$
and $\Hat h\colon A\times X\times[0,1]\to Y$. Being a closed subset of a metric space, $A$ is $P_0$-embedded
in $Z$. It follows that $A\times X$ is $P_0$-embedded in $Z\times X$, and therefore also in $U\times X$.
By Morita's Homotopy extension theorem, therefore, $\Hat g$ and $\Hat h$ induce a homotopy
$\Hat H\colon U\times X\times[0,1]\to Y$ extending both. The adjoint of $\Hat H$ is a continuous map
$H\colon U\times[0,1]\to Y^X$ which, on level $1$, is the desired extension of $f$. \qed

For general spaces $X$, the outline fails at two points, both of which are a consequence
of the failure of the exponential correspondence between continuous functions $Z\to Y^X$ and $Z\times X\to Y$
in case $X$ is not compactly generated. 

The need for the hypotheses on $X$ is the following: compactly generated Hausdorff is used for
a kind of exponential correspondence, and hemicompact is used to ensure that $Y^X$ is metrizable
(Fr\'echet) whenever $Y$ is metrizable (Fr\'echet).

For the rest of this section, let $X$ be a fixed compactly generated hemicompact space
with the `distinguished' sequence of compacta $\setof{K_i\vert\,i}$.

For the record, we cite the classical exponential correspondence theorem (see \cite{spanier}, Introduction, {\bf 8}).

\begin{prop}\label{cexp}
Let $X,Y,Z$ be topological spaces with $X$ locally compact Hausdorff (no separation properties are
assumed for $Y$ and $Z$). Let $Y^X$ be the space of continuous functions and let $f\colon Z\to Y^X$
be any function with set-theoretic adjoint $\Hat f\colon Z\times X\to Y$. Then $f$ is continuous if
and only if $\Hat f$ is continuous. This accounts for a bijection $(Y^X)^Z\leftrightarrow Y^{(X\times Z)}$. \qed
\end{prop}

\begin{defn}
For any space $Z$, let $\kappa(Z\times X)$ denote the topological space whose underlying set is $Z\times X$
and has its topology determined by the subsets $Z\times K_i$ (with the cartesian product topology). The identity
$\kappa(Z\times X)\to Z\times X$, where the latter has the cartesian product topology, is evidently continuous.
\end{defn}

The introduction of the topology $\kappa(Z\times X)$ is motivated by \lemref{exp}
whose proof is an easy consequence of \propref{cexp} together with the fact that for any space $Y$,
the space $Y^X$ is homeomorphic with the inverse limit $\lim_{i}Y^{K_i}$, as can easily be verified.

\begin{lem}\label{exp}
Let $Z$ and $Y$ be topological spaces (no separation axioms required) and let $X$ be a compactly generated hemicompact space.
Let $f\colon Z\to Y^X$ be a function with set-theoretic adjoint $\Hat f\colon X\times Z\to Y$. Then $f$ is continuous
if and only if $\Hat f\colon\kappa(Z\times X)\to Y$ is. This accounts for a bijection
$(Y^X)^Z\leftrightarrow Y^{\kappa(Z\times X)}$. \qed
\end{lem}

\begin{rem}
If $Z$ is compactly generated Hausdorff (as in our application) then in fact $\kappa(Z\times X)$ coincides with
the well known compactly generated refinement $\kkkk(Z\times X)$ of the cartesian product topology. If, in addition,
$Y$ is Hausdorff, the correspondence of \lemref{exp} is a homeomorphism. However, \lemref{exp} is all that we use,
and thus we do not need to recall properties of the functor $\kkkk$. Moreover, 
\propref{P} below is valid for arbitrary $P$-embeddings.
\end{rem}

\begin{lem}\label{cpt_prod}
Let $K$ be a compact Hausdorff space. Then $\kappa\big((Z\times K)\times X\big)$ is naturally
homeomorphic with $\kappa(Z\times X)\times K$.
\end{lem}
\begin{proof}
Applying the obvious bijection $Z\times K\times X\leftrightarrow Z\times X\times K$, we have to
show that the two topologies on $Z\times X\times K$ have the same continuous maps $Z\times X\times K$
(in fact in our application below we need exactly this fact). To this end, $f\colon \kappa(Z\times X\times K)\to Y$
is continuous if and only if the restrictions $f_i\colon Z\times K_i\times K\to Y$ are continuous which is if and only
if their adjoints $\Hat f_i\colon Z\times K_i\to Y^K$ are continuous. The latter is if and only if the map
$\Hat f\colon\kappa(Z\times X)\to Y^K$ is continuous and this in turn if and only if $f\colon\kappa(Z\times X)\times K\to Y$
is continuous which finishes the proof.
\end{proof}

Note that for a closed subset $A$ of $Z$ the topology $\kappa(A\times X)$
coincides with the topology that the set $A\times X$ inherits from $\kappa(A\times Z)$.
However, for arbitrary $A$ the two topologies might differ.

Let $(Z,A)$ be a topological pair (no separation properties assumed). Then $A$ is {\it $P$-embedded} in $Z$ if
continuous pseudo-metrics on $A$ extend to continuous pseudo-metrics on $Z$. Also, $A$ is a {\it zero set} in $Z$ if
there exists a continuous function $\phi\colon Z\to\RR$ with $A=\phi^{-1}(0)$. If $A$ is a $P$-embedded zero
set, it is called {\it $P_0$-embedded}.

For example, every closed subset of a metrizable space is $P_0$-embedded.

We need $P$-embeddings in the context of Morita's homotopy extension theorem:

\begin{thm}[Morita \cite{morita}]\label{morita_thm}
If $A$ is $P_0$-embedded in the topological space $Z$ then the pair
$(Z,A)$ has the homotopy extension property with respect to all ANR spaces. 
That is, if $Y$ is an ANR, if $g\colon Z\times\setof{0}\to Y$ and
$h\colon A\times[0,1]\to Y$ are continuous maps that agree pointwise on $A\times\setof{0}$,
then there exists a continuous map $H\colon Z\times[0,1]\to Y$ extending both $g$ and $h$.\qed
\end{thm}

\begin{prop}\label{P}
Let $A$ be $P$-embedded in $Z$ and let $X$ be a compactly generated hemicompact space. Then
the subset $A\times X$ is $P$-embedded in $\kappa(Z\times X)$.
\end{prop}

A result due to Al\`o and Sennott (see \cite{alo-sennott}, Theorem 1.2) shows that $A$ is $P$-embedded in $Z$
if and only if every continuous function from $A$ to a Fr\'echet space extends continuously over $Z$.
\propref{P} seems to be the right way of generalizing the equivalence (1)$\iff$(2) of Theorem 2.4 in \cite{alo-sennott}.

\begin{proof}
Let $E$ be a Fr\'echet space and let $f\colon A\times X\to E$ be a continuous map
where $A\times X$ is understood to inherit its topology from $\kappa(Z\times X)$.
Precomposing with the continuous identity $\kappa(A\times X)\to A\times X$ and using \lemref{exp},
we obtain a continuous map $A\to E^X$. As $E^X$ is also a Fr\'echet space (this is due to Arens,
see \cite{alo-sennott}, Proposition 2.2) and $A$ is $P$-embedded in $Z$, the function $\Hat f$
extends continuously to $\Hat F\colon Z\to E^X$. Reapplying \lemref{exp}, $\Hat F$ induces the
desired extension $F\colon\kappa(Z\times X)\to E$.
\end{proof}

\begin{proof}[Proof of \thmref{submain}]
Let $Z$ be metrizable and let $f\colon A\to Y^X$ be a continuous map defined on the closed subset
$A$ of $Z$. Let $g=g\colon U\times\setof{0}\to Y^X$ and $h$ be as in the above outline. Let
$h\sqcup g\colon A\times[0,1]\cup U\times\setof{0}\to Y^X$ denote the continuous union of the two,
with adjoint $\Hat h\sqcup\Hat g\colon\kappa\big((A\times[0,1]\cup U\times\setof{0})\times X\big)\to Y$.
As $A\times[0,1]\cup U\times\setof{0}$ is closed in $A\times U$, the map $\Hat h\sqcup\Hat g$ is
continuous with respect to the topology that $(A\times[0,1]\cup U\times\setof{0})\times X$ inherits
from $\kappa(U\times[0,1]\times X)$. Under the homeomorphism $\kappa(U\times[0,1]\times X)\approx
\kappa(U\times X)\times[0,1]$ of \lemref{cpt_prod}, the map $\Hat h\sqcup\Hat g$ corresponds to
$\Hat k\colon(A\times X)\times[0,1]\cup(U\times X)\times\setof{0}\to Y$.

Obviously, as $A$ is a zero set in $U$, the product $A\times X$ is a zero set in $U\times X$ with
respect to the cartesian product topology. {\it A fortiori}, $A\times X$ is a zero set in $\kappa(U\times X)$.
Hence, by \propref{P}, the set $A\times X$ is $P_0$-embedded in $\kappa(U\times X)$. \thmref{morita_thm}
yields an extension of
$\Hat k$ to $K\colon\kappa(U\times X)\times[0,1]\to Y$. Reapplying \lemref{cpt_prod} and \lemref{exp},
$K$ induces a continuous function $k\colon U\times[0,1]\to Y^X$. Level $1$ of this homotopy
is a continuous extension of $f$ over the neighbourhood $U$. Therefore, $Y^X$ is an ANR.
\end{proof}

\begin{cor}\label{not_metric}
Let $C$ be a compact Hausdorff space and let $Y$ be an ANR for metric spaces. Then $Y^C$ is an ANR for metric spaces. \qed
\end{cor}

\corref{not_metric} was proven independently by Yamashita \cite{yamashita} but the author of this note has not seen
it elsewhere except with the additional requirement of metrizability of $C$. From the point of view of $P$-embeddings,
however, \corref{not_metric} encodes a long-known fact (see \cite{alo-sennott}, Theorem 3.3): if $A$ is $P$-embedded in $Z$
and $X$ is compact Hausdorff, then $A\times X$ is $P$-embedded in $Z\times X$.


\begin{thebibliography}{9}

\bibitem{alo-sennott}
	Richard A.~Al\`{o}, Linnea Sennott,
	\booktitle{Collectionwise normality and the extension of functions on product spaces.}
	Fund.~Math.~ {\bf 76} (1972), no.~3, 231--243.

\bibitem{cauty_ANR} Robert Cauty,
	\booktitle{Une caract\'erisation des r\'etractes absolus de voisinage.}
	Fund.~Math.~ {\bf 144} (1994), no.~1, 11--22.

\bibitem{tdieck} Tammo tom Dieck,
	\booktitle{Partitions of unity in homotopy theory.}
	Compositio Math.~ {\bf 23} (1971), 159--167.

\bibitem{mardesic-segal} Sibe Marde\v{s}i\'{c}, Jack Segal,
	\booktitle{Shape theory. The inverse system approach.}
	North-Holland Mathematical Library, 26. North-Holland Publishing Co., Amsterdam-New York, 1982.

\bibitem{milnor} John Milnor,
	\booktitle{On spaces having the homotopy type of ${\rm CW}$-complex.}
	Trans.~Amer.~Math.~Soc.~ {\bf 90} (1959), 272--280. 

\bibitem{morita} Kiiti Morita,
	\booktitle{On generalizations of Borsuk's homotopy extension theorem.}
	Fund.~Math.~ {\bf 88} (1975), 1--6.

\bibitem{schon} Rolf Sch\"{o}n,
	\booktitle{Fibrations over a CWh-base.}
	Proc.~Amer.~Math.~Soc.~ {\bf 62} (1976), no. 1, 165--166 (1977).

\bibitem{smrekar} Jaka Smrekar,
	\booktitle{Compact open topology and CW homotopy type.}
	Topology and its applications {\bf 130} (2003), 291--304.

\bibitem{smr_CRM} Jaka Smrekar,
	\booktitle{CW type of inverse limits and function spaces.}
	Centre de Recerca Matem\`{a}tica, Bellaterra, Preprint no.~{\bf 548}, July 2003.
	Available electronically at \texttt{http://www.crm.es/Publications/03/pr548.pdf}
	and, in slightly revised version at \texttt{arXiv:0708.2838}

\bibitem{spanier} Edwin H. Spanier, \booktitle{Algebraic topology.}
	Corrected reprint of the 1966 original. Springer-Verlag, New York, 19??. 

\bibitem{stasheff} James Stasheff,
	\booktitle{A classification theorem for fibre spaces.}
	Topology {\bf 2} (1963), 239--246.

\bibitem{yamashita} Atsushi Yamashita,
	\booktitle{ANR mapping spaces with noncompact domains.}
	Preprint.


\end{thebibliography}
\end{document}